\documentclass{amsart}%
\usepackage{amsmath}
\usepackage{amsthm}
\usepackage{amssymb}
\usepackage{amsfonts}
\usepackage{graphicx}
\usepackage{url}%
\numberwithin{equation}{section}
\newtheorem{theorem}{Theorem}
\newtheorem{lemma}{Lemma}
\newtheorem{corollary}{Corollary}

\begin{document}
\title[On the Convergence of the Ensemble Kalman Filter]{On the Convergence \\ of the Ensemble Kalman Filter}
\author{Jan Mandel,
Loren Cobb,
Jonathan D. Beezley}
\begin{abstract}
Convergence of the ensemble Kalman filter in the limit for large ensembles to
the Kalman filter is proved. In each step of the filter, convergence of the
ensemble sample covariance follows from a weak law of large numbers for
exchangeable random variables, the continuous mapping theorem gives convergence 
in probability of the ensemble members, and $L^{p}$ bounds on the ensemble then give $L^{p}%
$ convergence.
\end{abstract}
\maketitle

\begin{center}{January 2009, revised May 2011}\end{center}
\bigskip

%\onehalfspacing

\section{Introduction}

Data assimilation uses statistical estimation to update the state of a running
model based on new data. Data assimilation is of great importance and widely
used in many disciplines including numerical weather prediction
\cite{Kalnay-AMD-2003}, ocean modeling \cite{Evensen-2007-DAE}, remote sensing
\cite{Vodacek-2008-RSD}, and image reconstruction \cite{Frazin-2005-TRN}. In
these applications, the dimension of the state is very high, often millions
and more, because the state consists of the values of a simulation on a
computational grid in a spatial domain. Consequently, the classical Kalman
filter (KF), which requires maintaining the state covariance matrix, is no
longer feasible.

One of the most successful recent data assimilation methods for
high-dimensional problems is the ensemble Kalman filter (EnKF). EnKF is a
Monte Carlo approximation of the KF, with the covariance in the KF\ replaced
by the sample covariance computed from an ensemble of realizations. Because
the EnKF does not need to maintain the state covariance matrix, it can be
implemented efficiently for high-dimensional problems. Although the EnKF
formulas rely on the assumption that the distribution of the state and the
data likelihood are normal, the ensemble can robustly describe an arbitrary
state probability distribution. Thus, in spite of errors such as smearing of
the state distribution towards normality \cite{Mandel-2009-EKP}, the EnKF is
often used for nonlinear systems.

One of the reasons for the popularity of the EnKF in applications is that the
convergence of EnKF with the ensemble size tends to be quite fast and
reasonably small ensembles (typically $25$ to $100$) are usually sufficient
\cite{Evensen-2007-DAE}. Convergence of the EnKF can be further accelerated by
localization, such as covariance tapering \cite{Furrer-2007-EHP}, which
improves the accuracy of the sample covariance. The EnKF
converges rapidly in these applications because the state vectors are not
arbitrary; rather, they are discretizations of smooth functions on a spatial
domain, and so they are the states of an infinitely dimensional dynamical system. One
explanation is that the state moves along a low-dimensional attractor. Indeed,
in weather simulations, the EnKF performance can be further improved by a
carefully chosen initial ensemble, which approximately covers the attractor well
\cite{Kalnay-AMD-2003}. Another explanation is that a smooth random field can
be well approximated by a linear combination of a small number of smooth
functions with random coefficients, such as a truncated random Fourier series
or Karhunen-Lo\`{e}ve expansion. Indeed, if the state is not smooth enough,
the convergence of the EnKF deteriorates \cite{Beezley-2009-DAC} and large
ensembles would be needed for acceptable accuracy.

A large body of literature on the EnKF and variants exists, but rigorous
probabilistic analysis is lacking. It is commonly assumed that the ensemble is
a sample (that is, i.i.d.) and that it is normally distributed. Although the
resulting analyses played an important role in the development of EnKF, both
assumptions are false. The ensemble covariance is computed from all ensemble
members together, thus introducing dependence, and the EnKF formula is a
nonlinear function of the ensemble, thus destroying the normality of the
ensemble distribution.

For example, the analysis in \cite{Burgers-1998-ASE} is based on the
comparison of the covariance of the analysis ensemble and the covariance of
the filtering distribution. The paper \cite{Furrer-2007-EHP} notes that if the
ensemble sample covariance is a consistent estimator, then Slutsky's theorem
yields the convergence in probability of the gain matrix. The paper
\cite{Li-2008-NPE} studies the interplay of numerical and stochastic errors.
All of these analyses assume that the ensemble covariance converges in some
sense in the limit for large ensembles, but a rigorous justification has not
yet become available.

This paper provides a rigorous proof that the EnKF converges to the KF in the
limit for large ensembles and for normal state probability distributions and
normal data likelihoods. The present analysis does not assume that the ensemble
members are independent or normally distributed. The ensemble members are
shown to be exchangeable random variables bounded in all $L^{p}$,
$p\in[1,\infty)$, which provides properties that replace independence and
normality. An argument using uniform integrability and the continuous mapping theorem is
then possible.

The result is valid for the EnKF version of Burgers, van Leeuven, and Evensen
\cite{Burgers-1998-ASE} in the case of constant state space dimension, a
linear model, normal data likelihood and initial state distributions, and
ensemble size going to infinity. This EnKF\ version involves randomization of
data. Efficient variants of EnKF without randomization exist
\cite{Anderson-2001-EAK,Tippett-2003-ESR}, but they are not the subject of
this paper.

Probabilistic analysis of the performance of the EnKF on nonlinear systems,
for non-normal state probability distributions, as well as analysis of the
speed of convergence of the EnKF to the KF and the dependence of the required
ensemble size on the state dimension, are outside of the scope of this paper
and left to future research. Some computational experiments and heuristic
explanations can be found in \cite{Beezley-2009-DAC}.

After the original preprint of this paper was completed \cite{Mandel-2009-CEK}, 
some related work became available. The
proof of EnKF convergence in \cite{Butala-2009-TID} has a gap; 
it assumes that certain covariances
derived from the ensemble exist, which is not guaranteed without an $L^{2}$
bound. The proof in \cite{LeGland-2009-LSA} is related and also uses a priori
$L^p$ bounds, but it appears
to be much longer and more complicated in order to obtain further analysis.

\section{Preliminaries}

The Euclidean norm of column vectors in $\mathbb{R}^{m}$, $m\geq1$, and the
induced matrix norm are denoted by $\Vert\cdot\Vert$, and $^{\mathrm{T}}$ is
the transpose. The stochastic $L^{p}$ norm of a random element $X$ is $\Vert
X\Vert_{p}=(E(\Vert X\Vert^{p}))^{1/p}$. The $j$-th entry of a vector $X$ is
$[X]_{j}$ and the $i,j$ entry of a matrix $Y\in\mathbb{R}^{m\times n}$ is
$[Y]_{ij}$. Convergence in probability is denoted by
$\xrightarrow{\mathrm{P}}$. We denote by
\[
X_{N}=[X_{Ni}]_{i=1}^{N}=[X_{N1},\ldots,X_{NN}],
\]
with various superscripts and for various $m\geq1$, an ensemble of $N$ random
elements in $\mathbb{R}^{m}$, called members. Thus, an ensemble is a random
$m\times N$ matrix with the ensemble members as columns. Given two ensembles
$X_{N}$ and $Y_{N}$, the stacked ensemble $[X_{N};Y_{N}]$ is defined as the
block random matrix%
\[
\lbrack X_{N};Y_{N}]=\left[
\begin{array}
[c]{c}%
X_{N}\\
Y_{N}%
\end{array}
\right]  =\left[  \left[
\begin{array}
[c]{c}%
X_{N1}\\
Y_{N1}%
\end{array}
\right]  ,\ldots,\left[
\begin{array}
[c]{c}%
X_{NN}\\
Y_{NN}%
\end{array}
\right]  \right]  =[X_{Ni};Y_{Ni}]_{i=1}^{N}.
\]
If all the members of $X_{N}$ are identically distributed, we write
$E(X_{N1})$ and $\operatorname*{Cov}(X_{N1})$ for their common mean vector and
covariance matrix. The ensemble sample mean and ensemble sample covariance
matrix are the random elements $\overline{X}_{N}=\frac{1}{N}\sum_{i=1}%
^{N}X_{Ni}$ and $C(X_{N})=\overline{X_{N}X_{N}^{\mathrm{T}}}-\overline{X}%
_{N}\overline{X}_{N}^{\mathrm{T}}$. \emph{All convergence is for
}$N\rightarrow\infty$\emph{.}

We will work with ensembles such that the joint distribution of the ensemble
$X_{N}$ is invariant under a permutation of the ensemble members. Such
ensemble is called \emph{exchangeable}. That is, an ensemble $X_{N}$, $N\geq
2$, is exchangeable if and only if $\Pr(X_{N}\in B)=\Pr(X_{N}\Pi\in B)$ for
every Borel set $B\subset\mathbb{R}^{m\times N}$ and every permutation matrix
$\Pi\in\mathbb{R}^{N\times N}$. The covariance between any two members of an
exchangeable ensemble is the same, $\operatorname*{Cov}(X_{Ni},X_{Nj}%
)=\operatorname*{Cov}(X_{N1},X_{N2})$, if $i\neq j$.

\begin{lemma}
\label{lem:exchangeable}Suppose $X_{N}$ and $D_{N}$ are exchangeable, the
random elements $X_{N}$ and $D_{N}$ are independent, and $Y_{Ni}=F(
X_{N},X_{Ni},D_{Ni}) $, $i=1,\ldots,N,$ where $F$ is measurable and
permutation invariant in the first argument, i.e. $F( X_{N}\Pi,X_{Ni},D_{Ni})
=F( X_{N},X_{Ni},D_{Ni}) $ for any permutation matrix $\Pi$. Then $Y_{N}$ is exchangeable.
\end{lemma}

\proof Write $Y_{N}=\mathbf{F}(X_{N},D_{N}),$ where
\[
\mathbf{F}(X_{N},D_{N})=[F(X_{N},X_{N1},D_{N1}),F(X_{N},X_{N2},D_{N2}%
),\ldots,F(X_{N},X_{NN},D_{NN})].
\]
Let $\Pi$ be a permutation matrix. Then $Y_{N}\Pi=\mathbf{F}(X_{N}\Pi,D_{N}%
\Pi)$. Because $X_{N}$ is exchangeable, the distributions of $X_{N}$ and
$X_{N}\Pi$ are identical. Similarly, the distributions of $D_{N}$ and
$D_{N}\Pi$ are identical. Since $X_{N}$ and $D_{N}$ are independent, the joint
distributions of $(X_{N},D_{N})$ and $(X_{N}\Pi,D_{N}\Pi)$ are identical.
Thus, for any Borel set $B\subset\mathbb{R}^{n\times N}$,
\begin{align*}
\Pr(Y_{N}\Pi\in B)  &  =E(1_{B}(Y_{N}\Pi))=E(1_{B}(\mathbf{F}(X_{N}\Pi
,D_{N}\Pi)))\\
&  =E(1_{B}(\mathbf{F}(X_{N},D_{N})))=\Pr(Y_{N}\in B),
\end{align*}
where $1_{B}$ stands for the characteristic function of $B$. Hence, $Y_{N}$ is
exchangeable. \endproof

We now prove a weak law of large numbers for nearly i.i.d. exchangeable ensembles.

\begin{lemma}
\label{lem:large-numbers}If for all $N$, $X_{N}$, $U_{N}$ are ensembles of
random variables, $[X_{N};U_{N}]$ is exchangeable, $\operatorname*{Cov}%
(U_{Ni},U_{Nj})=0$ for all $i\neq j$, $U_{N1}\in L^{2}$ is the same for all
$N$, and $X_{N1}\rightarrow U_{N1}$ in $L^{2}$, then $\overline{X}%
_{N}\xrightarrow{\mathrm{P}} E(U_{N1})$.
\end{lemma}

\proof Since $X_{N}$ is exchangeable, $\operatorname*{Cov}(X_{Ni}%
,X_{Nj})=\operatorname*{Cov}(X_{N1},X_{N2})$ for all $i,j=1,\ldots,N$, $i\neq
j$. Since $X_{N}-U_{N}$ is exchangeable, also $X_{N2}-U_{N2}\rightarrow0$ in
$L^{2}$. Then, using the identity $\operatorname*{Cov}(X,Y)=E(XY)-E(X)E(Y)$
and the Cauchy inequality for the $L^{2}$ inner product $E(XY)$, we have
\begin{align*}
\lefteqn{\left\vert \operatorname*{Cov}(X_{N1},X_{N2})-\operatorname
*{Cov}(U_{N1},U_{N2})\right\vert }\\
&  \qquad\leq2\Vert X_{N1}\Vert_{2}\Vert X_{N2}-U_{N2}\Vert_{2}+2\Vert
U_{N2}\Vert_{2}\Vert X_{N1}-U_{N1}\Vert_{2},
\end{align*}
so $\operatorname*{Cov}(X_{N1},X_{N2})\rightarrow 0$. By the same argument,
$\operatorname*{Var}(X_{N1})\rightarrow \operatorname*{Var}(U_{N1})<\infty$.
Now $E(\overline{X}_{N})=E(X_{N1})\rightarrow  E(U_{N1})$ from $X_{N1}%
-U_{N1}\rightarrow0$ in $L^{2}$, and
\begin{align*}
\operatorname*{Var}(\overline{X}_{N})  &  =\frac{1}{N^{2}}\sum_{i=1}%
^{N}\operatorname*{Var}(X_{Ni})+\sum_{i,j=1,j\neq j}^{N}\operatorname*{Cov}%
(X_{Ni},X_{Nj})\\
&  =\frac{1}{N}\operatorname*{Var}(X_{N1})+(1-\frac{1}{N})\operatorname*{Cov}%
(X_{N1},X_{N2})\rightarrow0,
\end{align*}
and the conclusion follows from the Chebyshev inequality. \endproof

The convergence of the ensemble sample covariance follows.

\begin{lemma}
\label{lem:conv-cov}If for all $N$, $X_{N}$, $U_{N}$ are ensembles of random
elements in $\mathbb{R}^{n}$, $\left[  X_{N};U_{N}\right]  $ is exchangeable,
$U_{N}$ are i.i.d., $U_{N1}\in L^{4}$ is the same for all $N$, and
$X_{N1}\rightarrow U_{N1}$ in $L^{4}$, then $\overline{X}_{N}\xrightarrow{\mathrm{P}}
E(U_{N1})$ and $C(X_{N})\xrightarrow{\mathrm{P}}\operatorname*{Cov}(U_{N1})$.
\end{lemma}

\proof From Lemma \ref{lem:large-numbers}, it follows that $[\overline{X}%
_{N}]_{j}\xrightarrow{\mathrm{P}}\lbrack E(U_{N1})]_{j}$ for each entry $j=1,\ldots,n$, so
$\overline{X}_{N}\xrightarrow{\mathrm{P}} E(U_{N1})$. Let $Y_{Ni}=X_{Ni}X_{Ni}^{\mathrm{T}%
}$, so that $C(X_{N})$ $=\overline{Y}_{N}-\overline{X}_{N}\overline{X}%
_{N}^{\mathrm{T}}$. Each entry of $[Y_{Ni}]_{j\ell}=[X_{Ni}]_{j}[X_{Ni}%
]_{\ell}$ satisfies the assumptions of Lemma \ref{lem:large-numbers}, so
$[Y_{Ni}]_{j\ell}\xrightarrow{\mathrm{P}} E([U_{N1}U_{N1}^{\mathrm{T}}]_{j\ell})$.
Convergence of the entries $[\overline{X}_{N}\overline{X}_{N}^{\mathrm{T}%
}]_{j\ell}=[\overline{X}_{N}]_{j}[\overline{X}_{N}]_{\ell}$ to $E([U_{N1}%
]_{j\ell})E([U_{N1}^{\mathrm{T}}]_{j\ell})$ follows from the already proved
convergence of $\overline{X}_{N}$ and the continuous mapping theorem \cite[p.
7]{vanderVaart-2000-AS}. Applying the continuous mapping theorem again, we get $C(X_{N}%
)\xrightarrow{\mathrm{P}}\operatorname*{Cov}(U_{N1})$. \endproof

\section{Formulation of the EnKF}

Consider an initial state given as the random variable $U^{(0)}$. In step $k$,
the state $U^{(k-1)}$ is advanced in time by applying the model $M^{(k)}$ to
obtain $U^{(k),f}=M^{(k)}(U^{(k-1)})$, called the prior or the forecast, with
probability density function (pdf) $p_{U^{(k),f}}$. The data in step $k$ are
given as measurements $d^{(k)}$ with a known error distribution, and expressed
as the data likelihood $p(d^{(k)}|u)$. The new state $U^{(k)}$ conditional on
the data, called the posterior or the analysis, then has the density
$p_{U^{(k)}}$ given by the Bayes theorem,
\[
p_{U^{(k)}}(u)\propto p(d^{(k)}|u)p_{U^{(k),f}}(u),
\]
where $\propto$ means proportional. This is the discrete-time filtering
problem. The distribution of $U^{(k)}$ is called the filtering distribution.

Assume that $U^{(0)}\sim N(u^{(0)},Q^{(0)})$, the model is linear,
$M^{(k)}:u\mapsto A^{(k)}u+b^{(k)}$, and the data likelihood is normal
conditional on given state $u^{(k),f}$,
\[
p\left(  d^{(k)}|u^{(k),f}\right)  \propto\mathrm{e}^{-\frac{1}{2}\left(
H^{(k)}u^{(k),f}-d^{(k)}\right)  ^{\mathrm{T}}R^{(k)^{-1}}\left(
H^{(k)}u^{(k),f}-d^{(k)}\right)  },
\]
where $H^{(k)}$ is the given observation matrix and $R^{(k)}$ is the given
data error covariance. The data error is assumed to be independent of the
model state. Then the filtering distribution is normal, $U^{(k)}\sim
N(u^{(k)},Q^{(k)})$, and it satisfies the KF
recursions~\cite{Anderson-1979-OF}%
\begin{align}
u^{(k),f}  &  =E(U^{(k),f})=A^{(k)}u^{(k)}+b^{(k)},\quad Q^{(k),f}%
=\operatorname*{Cov}U^{(k),f}=A^{(k)^{\mathrm{T}}}Q^{(k)}A^{(k)}%
,\label{eq:kf-rec1}\\
u^{(k)}  &  =u^{(k),f}+K^{(k)}(d^{(k)}-H^{(k)}u^{(k),f}),\quad Q^{(k)}%
=(I-K^{(k)}H^{(k)})Q^{(k),f}, \label{eq:kf-rec2}%
\end{align}
where the Kalman gain matrix $K^{(k)}$ is given by%
\begin{equation}
K^{(k)}=Q^{(k),f}H^{(k)\mathrm{T}}(H^{(k)}Q^{(k),f}H^{(k)\mathrm{T}}%
+R^{(k)})^{-1}. \label{eq:kf-gain}%
\end{equation}

The EnKF is obtained by replacing the exact covariance $Q^{(k)}$ by the
ensemble sample covariance and adding noise to the data in order to avoid a
shrinking of the ensemble spread and to obtain the correct filtering
covariance \cite{Burgers-1998-ASE}, cf. Lemma \ref{lem:kf-ens} below.

Let $U_{i}^{(0)}\sim N(u^{(0)},Q^{(0)})$ and $D_{i}^{(k)}\sim N(d^{(k)}%
,R^{(k)})$ be independent for all $k,i\geq1$. Given $N$, choose the initial
ensemble and the perturbed data as the first $N$ terms of the respective
sequence, $U_{Ni}^{(0)}=U_{i}^{(0)}$, $i=1,\ldots,N$, $D_{Ni}^{(k)}%
=D_{i}^{(k)}$, $i=1,\ldots,N$, $k=1,2,\ldots$The ensembles produced by EnKF
are $X_{N}^{(0)}=U_{N}^{(0)}$ and%
\begin{align}
X_{Ni}^{(k),f}  &  =\mathbf{\ }M^{(k)}(X_{Ni}^{(k-1)}),\quad i=1,\ldots
,N.\label{eq:M}\\
X_{N}^{(k)}  &  =X_{N}^{(k),f}+K_{N}^{(k)}(D_{N}^{(k)}-H^{(k)}X_{N}^{(k),f}),
\label{eq:enkf-random-step}%
\end{align}
where $K_{N}^{(k)}$ is the ensemble sample gain matrix,%
\begin{equation}
K_{N}^{(k)}=Q_{N}^{(k),f}H^{(k)T}(H^{(k)}Q_{N}^{(k),f}H^{(k)T}+R^{(k)}%
)^{-1},\quad Q_{N}^{(k),f}=C(X_{N}^{(k),f}). \label{eq:ens-gain}%
\end{equation}

Our analysis of the EnKF is based on the observation that the ensembles
$X_{N}^{(k)}$ are a perturbation of auxiliary ensembles $U_{N}^{(k)}$. The
ensembles $U_{N}^{(k)}$ are obtained from the same initial ensemble by
applying the KF formulas to each ensemble member separately and using the same
corresponding member of perturbed data,%
\begin{align}
U_{Ni}^{(k),f}  &  =\mathbf{\ }M^{(k)}(U_{Ni}^{(k-1)}),\quad i=1,\ldots
,N,\label{eq:M-U}\\
U_{N}^{(k)}  &  =U_{N}^{(k),f}+K^{(k)}(D_{N}^{(k)}-H^{(k)}U_{N}^{(k),f}).
\label{eq:ens-kf-step}%
\end{align}

The auxiliary ensembles $U_{N}^{(k)}$ are introduced for theoretical purposes
only and they do not play any role in the EnKF algorithm. The next lemma shows
that $U_{N}^{(k)}$ is a sample from the filtering distribution.

\begin{lemma}
\label{lem:kf-ens} For all $k=1,2,\ldots$, $U_{N}^{(k)}$ is i.i.d. and
$U_{N1}^{(k)}\sim N(u^{(k)},Q^{(k)})$.
\end{lemma}

\proof The statement is true for $k=0$ by definition of $U_{N}^{(0)}$. Assume
that it is true for $k-1$ in place of $k$. The ensemble $U_{N}^{(k)}$ is
i.i.d. and normally distributed, because it is an image under a linear map of
the normally distributed i.i.d. ensemble with members $[U_{Ni}^{(k-1)}%
,D_{Ni}^{(k)}]$, $i=1,\ldots,N$. Further, $D_{N}^{(k)}$ and $U_{Ni}^{(k),f}$
are independent, so from \cite[eq.~(15) and (16)]{Burgers-1998-ASE},
$U_{N1}^{(k)}$ has the correct mean and covariance, which uniquely determines the
normal distribution of $U_{N1}^{(k)}$. \endproof

\section{Convergence analysis}

\begin{lemma}
\label{lem:ens-Lp-bound}There exist constants $c( k,p) $ for all $k$ and all
$p \in[1,\infty)$ such that $\Vert X_{Ni}^{( k) }\Vert_{p}\leq c( k,p) $ and
$\Vert K_{N}^{( k) }\Vert_{p}\leq c( k,p) $ for all $N$.
\end{lemma}

\proof For $k=0$, each $X_{Ni}^{(k)}$ is normal. Assume $\Vert X_{Ni}%
^{(k-1)}\Vert_{p}\leq c(k-1,p)$ for all $N$. Then%
\[
\Vert X_{Ni}^{(k),f}\Vert_{p}=\Vert A^{(k)}X_{Ni}^{(k-1)}+b^{(k)}\Vert_{p}%
\leq\Vert A^{(k)}\Vert\Vert X_{Ni}^{(k-1)}\Vert_{p}+\Vert b^{(k)}\Vert
\leq\operatorname*{const}(k,p).
\]
By Jensen's inequality, for any $X_{N}$,
\[
\Vert\frac{1}{N}\sum_{i=1}^{N}X_{Ni}\Vert_{p}\leq\frac{1}{N}\sum_{i=1}%
^{N}\Vert X_{Ni}\Vert_{p}.
\]
This gives $\Vert\overline{X}_{N}^{(k),f}\Vert_{p}\leq\operatorname*{const}%
(k,p)$ and
\begin{align*}
\Vert Q_{N}^{(k),f}\Vert_{p}  &  \leq\frac{1}{N}\Vert X_{N1}^{(k),f}%
X_{N1}^{(k),f\mathrm{T}}\Vert_{p}+\frac{1}{N^{2}}\Vert X_{N1}^{(k),f}\Vert
_{p}^{2}\\
&  \leq\frac{1}{N}\Vert X_{N1}^{(k),f}\Vert_{2p}^{2}+\frac{1}{N^{2}}\Vert
X_{N1}^{(k),f}\Vert_{p}^{2}\leq\operatorname*{const}(k,p),
\end{align*}
since from the Cauchy inequality,
\begin{equation}
\Vert WZ\Vert_{p}\leq E\left(  \left\Vert W\right\Vert ^{p}\left\Vert
Z\right\Vert ^{p}\right)  ^{\frac{1}{p}}\leq E(\left\Vert W\right\Vert
^{2p})^{\frac{1}{2p}}E(\left\Vert Z\right\Vert ^{2p})^{\frac{1}{2p}%
}=\left\Vert W\right\Vert _{2p}\left\Vert Z\right\Vert _{2p},
\label{eq:cauchy}%
\end{equation}
for any compatible random matrices $W$ and $Z$. Since $H^{(k)}Q_{N}%
^{(k),f}H^{(k)\mathrm{T}}$ is symmetric positive semidefinite and $R^{(k)}$ is
symmetric positive definite, it holds that%
\[
\Vert(H^{(k)}Q_{N}^{(k),f}H^{(k)\mathrm{T}}+R^{(k)})^{-1}\Vert\leq
\Vert(R^{(k)})^{-1}\Vert\leq\operatorname*{const}(k),
\]
which, together with the bound on $\Vert Q_{N}^{(k),f}\Vert_{p}$, gives
\[
\Vert K_{N}^{(k)}\Vert_{p}\leq\Vert Q_{N}^{(k)}\Vert_{p}\operatorname*{const}%
(k)\leq\operatorname*{const}(k,p).
\]
Finally, we obtain the desired bound
\begin{align*}
\Vert X_{Ni}^{(k)}\Vert_{p}  &  \leq\Vert X_{Ni}^{(k),f}\Vert_{p}+\Vert
K_{N}^{(k)}D_{Ni}^{(k)}\Vert_{p}+\Vert K_{N}^{(k)}H^{(k)}X_{Ni}^{(k),f}%
\Vert_{p}\\
&  \leq\operatorname*{const}(k,p)(\Vert X_{Ni}^{(k),f}\Vert_{p}+\Vert
K_{N}^{(k)}\Vert_{p}+\Vert K_{N}^{(k)}\Vert_{2p}\Vert X_{Ni}^{(k),f}\Vert
_{2p})\leq c(k,p),
\end{align*}
using again (\ref{eq:cauchy}). \endproof

\begin{theorem}
\label{thm:conv-rand} For all $k$, $[X_{N};U_{N}]$ is exchangeable and
$X_{Ni}^{(k)}\rightarrow U_{Ni}^{(k)}$ in $L^{p}$ for all $p \in[1,\infty)$,
where $U_{N}$ is i.i.d. with the filtering distribution.
\end{theorem}

\proof The ensembles $U_{N}^{(k)}$ are obtained by linear mapping of the
i.i.d. initial ensemble $U_{N}^{(0)}$, so they are i.i.d. For $k=1$, we have $X_{N}%
^{(0)}=U_{N}^{(0)}$, $[X_{N}^{(0)};U_{N}^{(0)}]$ is exchangeable, and
$X_{Ni}=U_{Ni}$. Suppose the statement holds for $k-1$ in place of $k$. The
ensemble members are given by a recursion of the form
\[
\lbrack X_{Ni}^{(k)};U_{Ni}^{(k)}]=F^{(k)}(C(X_{N}^{(k-1)}),[X_{Ni}%
^{(k-1)};U_{Ni}^{(k-1)}],D_{Ni}^{(k)}).
\]
The ensemble sample covariance matrix $C(X_{N}^{(k-1)})$ is invariant to a 
permutation of ensemble members, so
$[X_{N}^{(k)};U_{N}^{(k)}]$ is exchangeable by Lemma \ref{lem:exchangeable}.
Since $X_{N}^{(k),f}$ and $U_{N}^{(k),f}$ satisfy the assumptions of Lemma
\ref{lem:conv-cov}, it follows that $C(X_{N}^{(k),f})\xrightarrow{\mathrm{P}}\operatorname*{Cov}%
U_{N1}^{(k),f}$ and $K_{N}^{(k)}\xrightarrow{\mathrm{P}} K^{(k)}$.
Thus, comparing (\ref{eq:enkf-random-step})\ and (\ref{eq:ens-kf-step}),  we have
that $X_{Ni}^{(k)}\xrightarrow{\mathrm{P}} U_{Ni}^{(k)}$,
by the continuous mapping theorem. Let $p \in[1,\infty)$. Since the sequence $\{X_{Ni}%
^{(k)}\}_{N=1}^{\infty}$ is bounded in $L^{p}$ by Lemma \ref{lem:ens-Lp-bound}
and $X_{Ni}^{(k)}\xrightarrow{\mathrm{P}} U_{Ni}^{(k)}$, it follows that $X_{Ni}%
^{(k)}\rightarrow U_{Ni}^{(k)}$ in $L^{q}$ for all $1\le q<p$ by uniform integrability \cite[p.
338]{Billingsley-1995-PM}. \endproof

Using Lemma \ref{lem:conv-cov} and uniform integrability again, it follows
that the ensemble mean and covariance converge to the filtering mean and covariance.

\begin{corollary}
$\overline{X}_{N}^{(k)}\rightarrow u^{(k)}$ and $C(X_{N}^{(k)})\rightarrow
Q^{(k)}$ in $L^{p}$ for all $p \in[1,\infty)$, where $u^{\left(  k\right)  }$
and $Q^{\left(  k\right)  }$ are the mean and the covariance of the filtering distribution.
\end{corollary}

\section{Acknowledgements}

This work was partially supported by the National Science Foundation under
grants CNS-0719641 and AGS-0835579 and the National Institutes of Health under
grant 1 RC1 LM01641-01. We would like to thank Professor Richard
Bradley for an example that led us away from mixing coefficients, and Dr.
Thomas Bengtsson for bringing \cite{Furrer-2007-EHP} to our attention.

%\newpage

%{\small
\bibliographystyle{abbrv}
\bibliography{enkf_theory}
%}

\bigskip

\emph{Authors' addresses}: \emph{Jan Mandel (corresponding author), Loren
Cobb, Jonathan D. Beezley}, Center for Computational Mathematics and
Department of Mathematical and Statistical Sciences, University of Colorado
Denver, Denver, CO 80217-3364, U.S.A.\newline e-mail:

\texttt{Jan.Mandel@ucdenver.edu}

\texttt{Loren.Cobb@ucdenver.edu}

\texttt{Jonathan.Beezley@ucdenver.edu}

\end{document}